\documentclass[a4paper, 11pt, reqno]{amsart}
\usepackage{a4wide}
\usepackage{amssymb, amsmath, amsthm}

\usepackage{hyperref}

\allowdisplaybreaks

\newtheorem{de}{Definition}[section]
\newtheorem{lem}[de]{Lemma}
\newtheorem{prop}[de]{Proposition}
\newtheorem{cor}[de]{Corollary}
\newtheorem{thm}[de]{Theorem}

\newtheorem{rem}[de]{Remark}

\numberwithin{equation}{section}

\newenvironment{prf}{\par\noindent{\bf Proof.\ }}{\hfill\rule{1.5ex}{1.5ex}\vspace{0.3cm}}

\newcommand{\bprf}{\begin{prf}}
\newcommand{\eprf}{\end{prf}}

\usepackage{amssymb}
\usepackage{dsfont}
\usepackage[T1]{fontenc}

\newcommand{\R}{\mathbb{R}}

\newcommand{\C}{\mathbb{C}}
\newcommand{\N}{\mathbb{N}}

\newcommand{\E}{\mathbb{E}}
\newcommand{\F}{\mathbb{F}}

\newcommand{\D}{\,\textrm{d}}

\newcommand{\trm}{\textrm}

\newcommand{\bea}{\begin{eqnarray*}}
\newcommand{\eea}{\end{eqnarray*}}

\newcommand{\beq}{\begin{equation}}
\newcommand{\eeq}{\end{equation}}

\newcommand{\tr}{\text{tr}}

\newcommand{\ra}{\rightarrow}

\newcommand{\hra}{\hookrightarrow}

\newcommand{\bc}{\begin{mathcal}}
\newcommand{\ec}{\end{mathcal}}
\newcommand{\calF}{\bc F\ec}

\newcommand{\calM}{\bc M\ec}
\newcommand{\calE}{\bc E\ec}
\newcommand{\calL}{\bc L\ec}

\newcommand{\calN}{\bc N\ec}
\newcommand{\calA}{\bc A\ec}
\newcommand{\calS}{\bc S\ec}

\newcommand{\calB}{\bc B\ec}
\newcommand{\calD}{\bc D\ec}

\newcommand{\dv}{\textrm{div}}

\newcommand{\ol}{\overline}

\newcommand{\wt}{\widetilde}

\newcommand{\rea}{\textrm{Re}\,}

\newcommand{\id}{\textrm{id}}

\newcommand{\eps}{\varepsilon}

\def\typeout#1{\message{^^J}\message{#1}\message{^^J}}
\newif\ifSRCOK \SRCOKtrue
\newcount\PAGETOP
\newcount\LASTLINE
\global\PAGETOP=1
\global\LASTLINE=-1
\def\EJECT{\SRC\eject}
\def\WinEdt#1{\typeout{:#1}}
\gdef\MainFile{\jobname.tex}
\gdef\CurrentInput{\MainFile}
\newcount\INPSP
\global\INPSP=0
\def\SRC{\ifSRCOK%
  \ifnum\inputlineno>\LASTLINE%
    \ifnum\LASTLINE<0%
      \global\PAGETOP=\inputlineno%
    \fi%
    \global\LASTLINE=\inputlineno%
    \ifnum\INPSP=0%
      \ifnum\inputlineno>\PAGETOP%
        
      \fi%
    \else%
      
    \fi%
  \fi%
\fi}
\def\PUSH#1{%
\SRC%
\ifnum\INPSP=0 \global\let\INPSTACKA=\CurrentInput \else%
\ifnum\INPSP=1 \global\let\INPSTACKB=\CurrentInput \else%
\ifnum\INPSP=2 \global\let\INPSTACKC=\CurrentInput \else%
\ifnum\INPSP=3 \global\let\INPSTACKD=\CurrentInput \else%
\ifnum\INPSP=4 \global\let\INPSTACKE=\CurrentInput \else%
\ifnum\INPSP=5 \global\let\INPSTACKF=\CurrentInput \else%
               \global\let\INPSTACKX=\CurrentInput \fi\fi\fi\fi\fi\fi%
\gdef\CurrentInput{#1}%
\WinEdt{<+ \CurrentInput}%
\global\LASTLINE=0%
\ifSRCOK\fi%
\global\advance\INPSP by 1}
\def\POP{%
\ifnum\INPSP>0 \global\advance\INPSP by -1  \fi%
\ifnum\INPSP=0 \global\let\CurrentInput=\INPSTACKA \else%
\ifnum\INPSP=1 \global\let\CurrentInput=\INPSTACKB \else%
\ifnum\INPSP=2 \global\let\CurrentInput=\INPSTACKC \else%
\ifnum\INPSP=3 \global\let\CurrentInput=\INPSTACKD \else%
\ifnum\INPSP=4 \global\let\CurrentInput=\INPSTACKE \else%
\ifnum\INPSP=5 \global\let\CurrentInput=\INPSTACKF \else%
               \global\let\CurrentInput=\INPSTACKX \fi\fi\fi\fi\fi\fi%
\WinEdt{<-}%
\global\LASTLINE=\inputlineno%
\global\advance\LASTLINE by -1%
\SRC}
\def\INPUT#1{\relax}

\let\originalxxxeverypar\everypar
\newtoks\everypar
\originalxxxeverypar{\the\everypar\expandafter\SRC}
\everymath\expandafter{\the\everymath\expandafter\SRC}
\output\expandafter{\expandafter\SRCOKfalse\the\output}

\newif\ifSRCOK \SRCOKtrue
\newcount\PAGETOP
\newcount\LASTLINE
\global\PAGETOP=1
\global\LASTLINE=-1
\gdef\MainFile{\jobname.tex}
\gdef\CurrentInput{\MainFile}
\newcount\INPSP
\global\INPSP=0
\def\EJECT{\SRC\eject}
\def\WinEdt#1{\typeout{:#1}}
\def\SRC{\ifSRCOK%
  \ifnum\inputlineno>\LASTLINE%
    \ifnum\LASTLINE<0%
      \global\PAGETOP=\inputlineno%
    \fi%
    \global\LASTLINE=\inputlineno%
    \ifnum\INPSP=0%
      \ifnum\inputlineno>\PAGETOP%
      \fi%
    \else%
    \fi%
  \fi%
\fi}
\def\PUSH#1{%
\SRC%
\ifnum\INPSP=0 \global\let\INPSTACKA=\CurrentInput \else%
\ifnum\INPSP=1 \global\let\INPSTACKB=\CurrentInput \else%
\ifnum\INPSP=2 \global\let\INPSTACKC=\CurrentInput \else%
\ifnum\INPSP=3 \global\let\INPSTACKD=\CurrentInput \else%
\ifnum\INPSP=4 \global\let\INPSTACKE=\CurrentInput \else%
\ifnum\INPSP=5 \global\let\INPSTACKF=\CurrentInput \else%
               \global\let\INPSTACKX=\CurrentInput \fi\fi\fi\fi\fi\fi%
\gdef\CurrentInput{#1}%
\WinEdt{<+ \CurrentInput}%
\global\LASTLINE=0%
\ifSRCOK\fi%
\global\advance\INPSP by 1}
\def\POP{%
\ifnum\INPSP>0 \global\advance\INPSP by -1  \fi%
\ifnum\INPSP=0 \global\let\CurrentInput=\INPSTACKA \else%
\ifnum\INPSP=1 \global\let\CurrentInput=\INPSTACKB \else%
\ifnum\INPSP=2 \global\let\CurrentInput=\INPSTACKC \else%
\ifnum\INPSP=3 \global\let\CurrentInput=\INPSTACKD \else%
\ifnum\INPSP=4 \global\let\CurrentInput=\INPSTACKE \else%
\ifnum\INPSP=5 \global\let\CurrentInput=\INPSTACKF \else%
               \global\let\CurrentInput=\INPSTACKX \fi\fi\fi\fi\fi\fi%
\WinEdt{<-}%
\global\LASTLINE=\inputlineno%
\global\advance\LASTLINE by -1%
\SRC}
\def\INPUT#1{\relax}
\let\OldINCLUDE=\include
\def\include#1{
\EJECT%
\PUSH{#1.tex}%
\OldINCLUDE{#1}%
\POP}

\let\originalxxxeverypar\everypar
\newtoks\everypar
\originalxxxeverypar{\the\everypar\expandafter\SRC}
\everymath\expandafter{\the\everymath\expandafter\SRC}
\let\zzzxxxbibliography=\bibliography
\def\bibliography#1{\PUSH{\jobname.bbl}\zzzxxxbibliography{#1}\POP}
\output\expandafter{\expandafter\SRCOKfalse\the\output}

\begin{document}

\title{Global attractors in stronger norms for a class of parabolic systems with nonlinear boundary conditions}

\author{Martin Meyries}

\address{Department of Mathematics,
Karlsruhe Institute of Technology, 76128 Karlsruhe, Germany.}
\email{martin.meyries@kit.edu}

\thanks{This paper is part of a research project supported by the Deutsche Forschungsgemeinschaft (DFG)}

\keywords{Parabolic systems, nonlinear boundary conditions, global attractors, maximal regularity, temporal weights, gradient estimates.}

\subjclass[2000]{35B41, 35K60}

\begin{abstract}\noindent For a class of quasilinear parabolic systems with nonlinear Robin boundary conditions we construct a compact local solution semiflow in a nonlinear phase space of high regularity. We further show that a priori estimates in lower norms are sufficient for the existence of a global attractor in this phase space. The approach relies on maximal $L_p$-regularity with temporal weights for the linearized problem. An inherent smoothing effect due to the weights is employed for gradient estimates. In several applications we can improve the convergence to an attractor by one regularity level.
\end{abstract}

\maketitle

\section{Introduction}\label{asyintro}
In this article we investigate the long-time behaviour of solutions in strong norms for nondegenerate reaction-diffusion systems with nonlinear Robin boundary conditions. 
For the unknown $u(t,x)\in \R^N$, where $N \in \N$, we consider (using sum convention)
\begin{alignat}{3}
 u_t - \partial_i(a_{ij}(u)\partial_j u) & =    f(u) & \qquad & \trm{in } \Omega, & \qquad &  t>0,  \nonumber\\
a_{ij}(u)\nu_i\partial_j u   & =  g(u) &&\trm{on } \Gamma,  &&  t>0,\quad  \label{nl1} \\
u|_{t=0} & =  u_0 &&  \trm{in } \Omega. && \nonumber 
\end{alignat} 
Here $\Omega\subset \R^n$ is a bounded domain with smooth boundary $\Gamma = \partial\Omega$, $n\geq 2$, and $\nu= (\nu_1,...\nu_n)$ denotes the outer normal field on $\Gamma$. We assume separated divergence form, i.e., 
$$a_{ij}(u) = a(u) \,\alpha_{ij}\in \R^{N\times N}, \qquad i,j\in \{1,...,n\},$$
where $a:\R^N \ra \R^{N\times N}$ and where the $\alpha_{ij}\in \R$ are constants, and further that
\beq\label{nl3}
\left.\begin{array}{c}
(\alpha_{ij})_{i,j=1,...,n}\text{ is symmetric and positive definite;}\\
 \text{for all $\zeta\in \R^N$ the spectrum of }a(\zeta) \text{ is contained in } \C_+= \{\rea z >0\}.
\end{array}\right\}
\eeq 
Systems with this structure are already considered in \cite{Ama90}. Moreover, $a$ and the reaction terms  $f,g:\R^N\ra \R^N$ are assumed to be smooth. Parabolic systems of type (\ref{nl1}) model a variety of phenomena in the sciences. An appropriate choice of $(a_{ij})$ yields, for instance, a system of heat equations, the Keller-Segel model for chemotaxis or a cross-diffusion population model (see Section 5). Observe that (\ref{nl3}) allows to rewrite the boundary condition in (\ref{nl1}) into the form $\alpha_{ij}\nu_i\partial_j u    =  \big(a^{-1}g\big)(u)$, which yields homogeneous Neumann conditions for $\alpha_{ij} = \delta_{ij}$ (the Kronecker symbol) and $g = 0$.

One often describes the long-time behaviour of solutions, as $t\ra +\infty$, in terms of a global attractor. Roughly speaking, a global attractor $\calA$ of the solution semiflow is a compact flow-invariant subset of the underlying phase space that attracts all bounded sets uniformly as $t\ra +\infty$; see  \cite{CD00} and \cite{Lad91} for the general theory. If $\calA$ is finite-dimensional then the complexity of the global dynamics of (\ref{nl1}) may essentially be reduced by restricting the semiflow to $\calA$. But although the solutions on $\calA$ may be smooth, it attracts only with respect to the metric of the phase space where the semiflow acts. It is therefore desireable to have an attractor in a phase space with metric  as strong as possible, say that of a Slobodetskii space $W_p^s$ with large $p$ and $s$ close to $2$. This can be useful, for instance, to improve error estimates for numerical algorithms when assuming in a quasi-stationary approximation that parts of a system of partial differential equations are on a fast time scale. 

Well-posedness, regularity and criteria for global existence of (\ref{nl1}) are  well understood. The general theory in \cite{Ama93} on quasilinear systems with nonlinear boundary conditions yields a solution semiflow in the phase spaces $W_p^s$ with sufficiently large $p$ and $s$ close to $1$. This approach, based on weak solutions, fits well to a priori estimates typically obtained in the applications, but in the end yields attractivity of $\calA$  with respect to a $C^\alpha$-norm, $\alpha\in (0,1)$. In particular, the long-time behaviour of the spatial gradient is not determined by $\calA$ in the sup-norm. Attractors in similiar norms are obtained in \cite{COPR97}, where the case of semilinear problems with nonlinear boundary conditions is treated. Things are simpler in the semilinear case with linear boundary conditions. If $\calA$ exists with respect to some $W_p^s$-norm with $s\in (0,2)$ then it is a consequence of the variation of constants formula that $\calA$ is attractive with respect to all $W_p^s$-norms, $s\in (0,2)$; see, e.g., Section 4.3 of \cite{CD00}.

The main point in this paper is thus to consider attractors for (\ref{nl1}) in a stronger norm. This leads to two difficulties.  First, if the norm is sufficiently strong then the  boundary conditions must hold in a trace sense, and thus nonlinear boundary conditions lead to a nonlinear phase space. Second,  compactness of the flow and an absorbant set are typically required for the existence of an attractor. To obtain this in a phase space of high regularity, a priori estimates in strong norms must be found for (\ref{nl1}), which is a rather delicate issue in many applications, especially when dealing with systems. 

Concerning the nonlinear phase space, in \cite{LPS06} a local semiflow in 
$$\calM_p^{2-2/p} = \big\{ u_0\in W_p^{2-2/p}(\Omega,\R^N)\;:\; a_{ij}(u_0)\nu_i\partial_j u_0    =  g(u_0) \big\}, \qquad p\in (n+2,\infty),$$
has been constructed and the local dynamics around an equilibrium have been discussed. These results rely on maximal $L_p$-regularity for the linearized problem with inhomogeneous boundary conditions, as treated in \cite{DHP07}. In the present article we can overcome the difficulty of a priori estimates in strong norms and consider the global dynamics of (\ref{nl1}) in $\calM_p^{2-2/p}$ in terms of attractors. Our main results may be summarized as follows. For a precise definition of a local semiflow and a global attractor we refer to the Sections 3 and 4, respectively.

\begin{thm}\label{thm}\textsl{Suppose that (\ref{nl3}) is valid, that the nonlinearities  are smooth, and take $p>n+2$. Then (\ref{nl1}) generates a compact local semiflow of solutions in $\calM_p^{2-2/p}$, such that for $u_0\in \calM_p^{2-2/p}$ the corresponding maximal solution $u(\cdot,u_0)$ belongs to 
$$W_{p}^1\big (0,\tau; L_p(\Omega,\R^N) \big) \cap L_p\big (0,\tau; W_p^2(\Omega;\R^N)\big),$$
for all $\tau\in \big(0,t^+(u_0)\big)$, with the maximal existence time $t^+(u_0)>0$. If the semiflow has an absorbant ball with respect to a H\"older norm, i.e., there are $\alpha, R>0$ such that 
\beq\label{n1}
\limsup_{t\ra t^+(u_0)} |u(t,u_0)|_{C^\alpha(\ol{\Omega},\R^N)} \leq R
\eeq 
for all $u_0\in \calM_p^{2-2/p}$ then (\ref{nl1}) has a global attractor in $\calM_p^{2-2/p}$. In the semilinear case, i.e., if $(a_{ij})$ is independent of $u$, it is sufficient to have an absorbant ball with respect to $W_q^\sigma\cap C(\ol{\Omega},\R^N)$ for some $\sigma>0$ and $q\in (1,\infty)$.}
\end{thm}

 We emphasize that our semiflow for (\ref{nl1}) is in any case compact in $\calM_p^{2-2/p}$, and that our method allows to show this also for more general problems (see Remark \ref{remark}). Moreover, an a priori estimate in a H\"older norm is sufficient for the existence of an attractor in this phase space of high regularity. In particular, if an absorbant set is known in a $W_p^1$-norm with $p>n$ then the theorem applies. In special situations, like single equations or triangular cross-diffusion systems, we can lower the strength of the metric for the absorbing ball even more by employing De Giorgi-Nash-Moser theory (see Section 4). Due to  $\calM_p^{2-2/p}\hra C^{1+\beta}$ for given $\beta\in (0,1)$ if $p$ is large enough, the theorem can give long-time control of the attractor over the gradient of solutions in a H\"older norm. In Section 5 we consider applications to a system of heat equations, a chemotaxis model and a cross-diffusion population model and improve the known convergence to an attractor. Here we use the results of \cite{Dun97}, \cite{Dun00} and \cite{KD07}.

Our results rely on a maximal $L_p$-regularity approach with temporal weights for linear parabolic problems with inhomogeneous boundary conditions, developed in \cite{MS11b}. The temporal regularity in this approach is based on the spaces
$$L_{p,\mu}(\R_+; E) = \big\{ u: \R_+\ra E\;:\; [t\mapsto t^{1-\mu} u(t)]\in L_p(\R_+; E)\big\}, \qquad \mu\in (1/p,1], \quad p\in (1,\infty),$$
where $E$ is a Banach space. The fact that the weight $t^{1-\mu}$ vanishes at $t=0$, and only there, allows for flexibility in the initial regularity in the maximal regularity approach for the linearization of (\ref{nl1}), and builds an inherent smoothing effect into the solutions. Well-posedness in a scale of compactly embedded nonlinear phase spaces is obtained by linearization, a detailed study of nonlinear superposition operators on weighted anisotropic spaces and the contraction principle. The smoothing effect due to the weights is used to show compactness properties of the semiflow in $\calM_p^{2-2/p}$ and to establish a gradient estimate. 

We briefly sketch the idea how this works in a linear situation. Let $-A$ be the generator of an exponentially stable analytic $C_0$-semigroup on $E$ with domain $D(A)$. Basic interpolation arguments show that $e^{-\cdot A} u_0$ belongs to $\E_{1,\mu}(\R_+):= W_{p,\mu}^1(\R_+; E)\cap L_{p,\mu}(\R_+; D(A))$ if and only if $u_0$ belongs to the real interpolation space  $(E,D(A))_{\mu-1/p,p}$, and in this case $|e^{-\cdot A} u_0|_{\E_{1,\mu}(\R_+)}\lesssim |u_0|_{(E,D(A))_{\mu-1/p,p}}$. On the other hand, since the weight does not vanish for positive times, the temporal trace $\tr_{t=\tau}$ is for $\tau>0$ continuous from $\E_{1,\mu}(\R_+)$ to $(E,D(A))_{1-1/p,p}$. Combining these facts and assuming that $E=L_p(\R^n)$ and $D(A)= W_p^2(\R^n)$, we obtain
$$|e^{-\tau A}u_0|_{W_p^{2-2/p}(\R^n)} \lesssim |u_0|_{W_p^{2\mu-2/p}(\R^n)}.$$
In this way one can control the solution of a linear evolution equation in a strong norm by its initial value in a lower norm: observe that $W_p^{2\mu-2/p}$ tends to $L_p$ as $\mu$ tends to $1/p$.

The importance of the $L_{p,\mu}$-spaces in the context of maximal regularity for linear problems has first been observed in \cite{PS04} and has  been used in \cite{KPW10} to show compactness properties of the semiflow for quasilinear problems with linear boundary conditions. Here we use it also for uniform gradient estimates: the extension of the above reasoning in the linear case via weighted maximal regularity to nonlinear problems (see Lemma \ref{nl31}) should be seen as the main technical contribution of the present paper.

This article is organized as follows. In Section 2 we introduce weighted anisotropic function spaces and study the properties of superposition operators associated to (\ref{nl1}) on them. In Section 3 we construct the compact local solution semiflow using linearization and the contraction principle. Gradient estimates are shown in Section 4, and applications to attractors of concrete models from the sciences are given in Section 5.

\textbf{Notation.} Although dealing with systems, we often write $L_p(\Omega) = L_p(\Omega,\R^N)$ and similiarly for other function spaces. We further write $a\lesssim b$ for some quantities $a,b$ if there is a generic positive constant $C$ with $a\leq Cb$. If $X,Y$ are Banach spaces we denote by $\calB(X,Y)$ the space of bounded linear operators between them, with $\calB(X) := \calB(X,X)$.

\section{Weighted function spaces and superposition operators}
For a Banach space $E$, a finite or infinite interval $J=(0,T)$, $p\in (1,\infty)$ and $\mu\in (1/p,1]$ we work with the weighted spaces
$$L_{p,\mu}(J;E) = \big\{u\;:\; t^{1-\mu}u \in L_p(J;E)\big\}, \qquad W_{p,\mu}^1(J;E) = \big\{u\;:\; u,u'\in L_{p,\mu}(J;E)\big\},$$
equipped with their canonical norms. We look for solutions of (\ref{nl1}) in the anisotropic space
$$\E_{1,\mu}(J) : = W_{p,\mu}^1\big(J; L_p(\Omega)\big)\cap L_{p,\mu}\big(J; W_p^{2}(\Omega)\big),$$
which suggests that the basic space for the domain equation in (\ref{nl1}) equals
$$\E_{0,\mu}(J): = L_{p,\mu}\big(J; L_p(\Omega)\big).$$ 
These function spaces are discussed in \cite{MS11a}. Denoting by $\tr_\Omega$ the spatial trace operator on $\Omega$, i.e., $\tr_\Omega u= u|_{\Gamma}$, Lemma 3.4 and Theorem 4.5 of \cite{MS11a} show that the Neumann boundary operator $\tr_\Omega \partial_j$ maps $\E_{1,\mu}(J)$ continuously into 
$$\F_\mu(J) := W_{p,\mu}^{1/2-1/2p}\big (J; L_p(\Gamma)\big) \cap L_{p,\mu}\big(J; W_p^{1-1/p}(\Gamma)\big).$$
Thus $\F_\mu(J)$ is the basic space for the boundary equation of (\ref{nl1}). Here $W_{p,\mu}^s(J;E)$ is for $s>0$ a weighted Slobodetskii space, which is defined by real interpolation between $L_{p,\mu}$ and $W_{p,\mu}^1$. An equivalent intrinsic norm for $W_{p,\mu}^s$ is given by $|u|_{L_{p,\mu}(J;E)} + [u]_{W_{p,\mu}^{s}(J; E)}$, where 
\beq\label{n3}
[u]_{W_{p,\mu}^{s}(J; E)}^p = \int_0^T \int_0^\tau  \frac{t^{p(1-\mu)}}{(\tau-t)^{1+s p}} |u(\tau) - u(t)|^p  \D t \D \tau.
\eeq 
Here the equivalence constants depend on $J=(0,T)$ and tend to infinity as $T$ tends to zero. This technical point becomes relevant  when working with short time intervals. We refer to \cite{MS11a} for more properties of these weighted spaces. Further, $W_p^{\kappa}(\Gamma)$ denotes for $\kappa>0$ a Slobodetskii space over the boundary $\Gamma$, which is defined by local charts; see e.g. Definition 3.6.1 in \cite{Tri94}. 
Theorem 4.2 of \cite{MS11a} gives the embedding
$$\E_{1,\mu}(J) \hra BU\!C\big(\ol{J}; B_{p,p}^{2(\mu-1/p)}(\Omega)\big).$$
Here $B_{p,p}^\kappa$ denotes a Besov space of order $\kappa>0$, that satisfies $W_p^\kappa = B_{p,p}^\kappa$ for $\kappa\notin \N_0$ (see again \cite{Tri94}). Therefore, by Sobolev's embeddings,
\beq\label{nl6}
\E_{1,\mu}(J)\hra BU\!C\big(\ol{J}; C^1(\ol{\Omega})\big) \qquad\text{if }\; 2(\mu-1/p) > 1+n/p.
\eeq
Similarly, it holds
$$\F_{\mu}(J) \hra BU\!C\big (\ol{J}; B_{p,p}^{2(\mu-1/p) -1-1/p}(\Gamma)\big) \qquad \text{if }\; 2(\mu-1/p) > 1+1/p,$$
so that we have
\beq\label{nl7}
\F_{\mu}(J) \hra BU\!C\big (\ol{J}\times \Gamma\big) \qquad \text{if }\; 2(\mu-1/p) > 1+n/p.
\eeq
Due to Lemma 4.3 of \cite{MS11a} there is a continuous right-inverse $S$ of the temporal trace $\tr_0: \E_{1,\mu}(\R_+) \ra B_{p,p}^{2(\mu-1/p)}(\Omega),$ i.e., $\tr_0 u = u|_{t=0}$. Observe that the relation $2(\mu-1/p) > 1+n/p$ for some $\mu\in (1/p,1]$ implies that $p>n+2$.

We also work with weighted spaces based on vanishing initial values and set for $s\in (0,1)$
$${}_0W_{p,\mu}^1(J;E) = \{u\in W_{p,\mu}^1(J;E): u(0) = 0\}, \qquad {}_0W_{p,\mu}^s(J;E) = \big(L_{p,\mu}(J;E),{}_0W_{p,\mu}^1(J;E)\big)_{s,p}.$$  By Proposition 2.10 of \cite{MS11a}, the temporal trace at $t=0$ is defined and continuous on $W_{p,\mu}^s$ if $s>1-\mu+1/p$, and it holds that
$${}_0W_{p,\mu}^s = W_{p,\mu}^s \;\text{ if } s <1-\mu+1/p,\qquad  {}_0W_{p,\mu}^s = \big\{ u\in W_{p,\mu}^s\,:\, u(0) = 0\big\}\;\text{ if } s >1-\mu+1/p.$$ 
In particular, ${}_0W_{p,\mu}^s$ is a closed subspace of $W_{p,\mu}^s$ for $s\neq 1-\mu+1/p,$ and the intrinsic norm (\ref{n3}) is also an equivalent norm.  Replacing the $W_{p,\mu}^s$-spaces in the definition of $\E_{1,\mu}$ and $\F_{\mu}$ by ${}_0W_{p,\mu}^s$-spaces, we denote the resulting spaces by ${}_0\E_{1,\mu}$ and ${}_0\F_{\mu}$, respectively. It is shown in \cite{MS11a} that if one restricts to ${}_0\E_{1,\mu}$ and ${}_0\F_{\mu}$ in the above embeddings, the embedding constants are independent of $T$. Moreover, by Lemma 2.5 of \cite{MS11a},  for $s\in [0,1]$ there is a continuous extension operator $\calE_J^0: {}_0W_{p,\mu}^s(J;E) \ra {}_0W_{p,\mu}^s(\R_+;E)$ whose norm is independent of the length of $J$. In a canonical way $\calE_J^0$ induces  an extension operator for ${}_0\E_{1,\mu}(J)$ and ${}_0\F_{\mu}(J)$ to the half-line.

Let us now study the properties of the nonlinear superposition operators occurring in (\ref{nl1}) in this weighted setting. We first consider the map $A$, defined by 
$$A(u) =  - \big(  \partial_i(a_{ij}(u)\partial_j u ) + f(u)\big), \qquad u\in \E_{1,\mu}(J).$$

\begin{lem}\label{nl19}\textsl{Let $J=(0,T)$ be finite, and let $p\in (n+2,\infty)$ and $\mu\in (1/p,1]$ be such that $2(\mu-1/p) > 1+n/p$. Then 
$A\in C^1\big(\E_{1,\mu}(J), \E_{0,\mu}(J)\big),$ and for $u\in \E_{1,\mu}(J)$ we have
$$-A'(u)h =  \partial_i(a_{ij}(u)\partial_j h + a_{ij}'(u)\partial_j uh ) + f'(u)h, \qquad h\in \E_{1,\mu}(J).$$
Moreover, let $T_0,R>0$ be given. Then there is a continuous function $\eps:[0,\infty)\ra [0,\infty)$ with $\eps(0) =0$ such that for $T\leq T_0$ it holds
$$
|A(u+h)-A(u) -A'(u)h|_{\E_{0,\mu}(J)}\leq \eps(|h|_{\E_{1,\mu}(J)})|h|_{\E_{1,\mu}(J)}
$$
for all $u, h \in \E_{1,\mu}(J)$ with 
\beq\label{nl41}
h|_{t=0} = 0, \qquad |u|_{C(\ol{J}; C^1(\ol{\Omega}))}, |u|_{\E_{1,\mu}(J)}, |h|_{\E_{1,\mu}(J)}\leq R.
\eeq}
\end{lem}
\bprf Using  $C\big(\ol{J}\times \ol{\Omega}\big)\hra \E_{0,\mu}$ and the embedding (\ref{nl6}), standard estimates show that for  $u,h\in \E_{1,\mu}(J)$ we have
\begin{align}
|A(u+h)&-A(u) - A'(u)h|_{\E_{0,\mu}(J)}  \nonumber \\
&\, \lesssim |f(u+h)-f(u) -f'(u)h|_{C(\ol{J}\times \ol{\Omega})} + |a_{ij}'(u)|_{C(\ol{J}; C^1(\ol{\Omega}))} |h|_{\E_{1,\mu}(J)}^2 \label{nl20}\\
&\, \qquad + |a_{ij}(u+h)-a_{ij}(u) -a_{ij}'(u)h|_{C(\ol{J}; C^1(\ol{\Omega}))}(|u|_{\E_{1,\mu}(J)} + |h|_{\E_{1,\mu}(J)}),\nonumber
\end{align}
where the maximum over single indices is understood. The differentiability of $f$ implies that
$$|f(u+h)-f(u) -f'(u)h|_{C(\ol{J}; C(\ol{\Omega}))} \leq \eps(|h|_{C(\ol{J}\times \ol{\Omega})} )|h|_{C(\ol{J}\times \ol{\Omega})} \leq \eps(|h|_{\E_{1,\mu}(J)} )|h|_{\E_{1,\mu}(J)}.
$$
In case (\ref{nl41}), the images of $u$ and $h$ are contained in a fixed compact subset of $\R^N$, which yields that $\eps$ is uniform in $T\leq T_0$ and $R$. The second  summand in (\ref{nl20}) may be estimated by $\eps(|h|_{\E_{1,\mu}(J)})|h|_{\E_{1,\mu}(J)}$, where $\eps$ is again uniform for (\ref{nl41}). For the third summand we have that the second factor is bounded, and it is uniformly bounded for (\ref{nl41}). Using (\ref{nl6}), the first factor there may be estimated in a standard way by $\eps(|h|_{\E_{1,\mu}(J)} )|h|_{\E_{1,\mu}(J)}$, again uniformly in $T\leq T_0$ and $R$ for (\ref{nl41}). This shows the differentiability of $A$ and the asserted uniformity of the linear approximation. Similiar considerations yield the continuity of the derivative $A'.$ \eprf

We next investigate the boundary nonlinearities. To this end we define for $q\in (1,\infty)$, $\mu\in (1/q,1]$ and $\kappa, \tau \in (0,1)$ the spaces
$$W_{q,\mu}^{\kappa, \tau }(J\times \Gamma) := W_{q,\mu}^\kappa\big( J; L_q(\Gamma)\big) \cap L_{q,\mu}\big( J; W_q^{\tau}(\Gamma)\big).$$
Of particular importance is the estimate (\ref{n101}) below which is useful for low values of $q$ and $\mu$.

\begin{lem}\label{jesaba1}\textsl{Let $J=(0,T)$ be finite and $g:\R^N \ra \R^N$ be smooth. Then for $q\in  (1,\infty)$, $\mu\in (1/q,1]$ and $\kappa,\tau \in (0,1)$ it holds
\beq\label{n101}
|g(u)|_{W_{q,\mu}^{\kappa, \tau }(J\times \Gamma)} \lesssim  \sup_{\zeta\in B_u}|g'(\zeta)| \, |u|_{W_{q,\mu}^{\kappa,\tau }(J\times \Gamma)} + |g(u)|_{C(\ol{J}\times \Gamma)}
\eeq
for all $u\in  W_{q,\mu}^{\kappa, \tau } \cap C(\ol{J}\times \Gamma)$, where $B_u$ is a ball with $u(\ol{J}\times \Gamma)\subset B_u$. Let now $p\in (n+2,\infty)$ and $\mu\in (1/p,1]$ satisfy $2(\mu-1/p)>1+n/p$. Then for the superposition operator $G$, given by $G(u) := g(\emph{\tr}_\Omega u)$, we have $G\in C^1\big (\E_{1,\mu}(J),\F_\mu(J)\big)$, with $G'(u) = g'(\emph{\tr}_\Omega u)\emph{\tr}_\Omega.$
Moreover, if $T_0,R>0$ are given, then there is a continuous function $\eps:[0,\infty)\ra [0,\infty)$ with $\eps(0)=0$ such that for $T\leq T_0$ it holds
\beq\label{n102}
|g(u+h)-g(u)-g'(u)h|_{ {}_0\F_\mu(J)} \leq \eps(|h|_{\E_{1,\mu}(J)}) |h|_{\E_{1,\mu}(J)}
\eeq
for all $u,h\in \E_{1,\mu}(J)$ satisfying
\beq\label{nl43}
h|_{t=0} = 0,\qquad  |u|_{C(\ol{J}; C^1(\Gamma))}, |u|_{\E_{1,\mu}(J)},\big|u|_{t=0}\big|_{B_{p,p}^{2(\mu-1/p)}(\Omega)},|h|_{\E_{1,\mu}(J)}\leq R.
\eeq}
\end{lem}

\bprf \textbf{(I)} To show (\ref{n101}) we take $u\in W_{q,\mu}^{\kappa, \tau}\cap C(\ol{J}\times \Gamma)$. Then it is clear that
\begin{align*}
|g(u)|_{L_{q,\mu}(J;L_q(\Gamma))} \lesssim |g(u)|_{C(\ol{J}\times \Gamma)}.
\end{align*}
To estimate $|g(u)|_{W_{q,\mu}^{\kappa}(J; L_q(\Gamma))}$ we use the intrinsic norm for $W_{q,\mu}^{\kappa}$ from (\ref{n3}), for which the mean value theorem immediately gives
$$[g(u)]_{W_{q,\mu}^{\kappa}(J; L_q(\Gamma))}^q \leq  \sup_{\zeta\in B_u}|g'(\zeta)|^q \;[u]_{W_{q,\mu}^{\kappa}(J; L_q(\Gamma))}^q.$$ 
Using a partition of unity and the intrinsic norm for $W_q^\tau$ given by Remark 4.4.1/2 of \cite{Tri94}, we obtain in the same way that
$$|g(u(t,\cdot))|_{W_q^{\tau}(\Gamma)} \lesssim \sup_{\zeta\in B_u}|g'(\zeta)| [u(t,\cdot)]_{W_q^{\tau}(\Gamma)} + |g(u)|_{C(\ol{J}\times \Gamma)}$$
for $t\in J$. Taking the $L_{q,\mu}$-norm yields (\ref{n101}).

\textbf{(II)} We next consider the map $G$. For $u\in\E_{1,\mu}(J)$ we have  from Theorems 4.2 and 4.5 of \cite{MS11a}, $2(\mu-1/p)> 1+n/p$ and $p>n$ that 
\beq\label{n2}
\tr_\Omega u\in W_{p,\mu}^{1-1/2p, 2-1/p}(J\times \Gamma) \hra C(\ol{J}\times \Gamma) \cap L_{p,\mu}\big (J; C^1(\Gamma)\big).
\eeq
In particular, it holds $g(\tr_\Omega u), g'(\tr_\Omega u)\in \F_\mu(J) \cap C(\ol{J}\times \Gamma)$ due to (\ref{n101}). Using (\ref{n2}) and the intrinsic norms from above for $\F_\mu(J)$, one easily obtains that $g'(\tr_\Omega u)\tr_\Omega \in \calB\big (\E_{1,\mu}(J),\F_\mu(J)\big ).$
Hence the maps $G$ and $G'$ are well-defined. To show the differentiability of $G$ at $u$, take $h\in \E_{1,\mu}(J)$. In the sequel we neglect the trace $\tr_\Omega$. It follows from standard arguments that there is $\eps:[0,\infty)\ra [0,\infty)$ with $\eps(0)=0$, which is uniform in $t\in J$ and $R$ for (\ref{nl43}), such that
\begin{align*}
 |g\big (u(t,\cdot)+h(t,\cdot)\big)-g\big(u(t,\cdot)\big)-g'\big(u(t,\cdot)\big)&\,h(t,\cdot)|_{C^1(\Gamma)} \leq \eps(|h(t,\cdot)|_{C^1(\Gamma)}) |h(t,\cdot)|_{C^1(\Gamma)}.
\end{align*}
Taking the $L_{p,\mu}$-norm, using the embeddings $C^1(\Gamma)\hra W_p^{1-1/p}(\Gamma)$ and (\ref{nl6}), we obtain 
\begin{align*}
 |g(u+h)-g(u)-g'(u)h|_{L_{p,\mu}(J;W_p^{1-1/p}(\Gamma))}&\, \leq \eps(|h|_{C(\ol{J}; C^1(\Gamma))}) |h|_{C(\ol{J}; C^1(\Gamma))}\\
&\, \lesssim \eps(|h|_{\E_{1,\mu}(J)}) |h|_{\E_{1,\mu}(J)}.
\end{align*}
Observe that these estimates are always uniform in $T\leq T_0$ and $R$ if (\ref{nl43}) holds. For the intrinsic seminorm of $W_{p,\mu}^{1/2-1/2p}\big (J; L_p(\Gamma)\big )$ we set 
$$\Xi(t,x):= g\big(u(t,x)+h(t,x)\big)-g\big(u(t,x)\big) - g'\big(u(t,x)\big)h(t,x),$$ 
and estimate with the mean value theorem
\begin{align}
[g(u+h)-g(u&) - g'(u)h]_{W_{p,\mu}^{1/2-1/2p}(J; L_p(\Gamma))}^p\nonumber \\
&= \int_0^T\int_0^s \int_\Gamma \frac{t^{p(1-\mu)}}{(s-t)^{1+ (1/2-1/2p)p}} |\Xi(s,x)-\Xi(t,x)|^p \D \sigma(x)\D t\D s\nonumber \\
&\, \leq \eps(|h|_{C(\ol{J}\times \Gamma)}^p)\big([h]_{W_{p,\mu}^{1/2-1/2p}(J; L_p(\Gamma))}^p + |h|_{C(\ol{J}\times \Gamma)}^p[u]_{W_{p,\mu}^{1/2-1/2p}(J; L_p(\Gamma))}^p \big)\nonumber\\
&\,\lesssim \eps( |h|_{\E_{1,\mu}(J)}^p) |h|_{\E_{1,\mu}(J)}^p. \label{nl45}
\end{align}
Therefore $G$ is differentiable. Similiar estimates yield the continuity of $G'$. 

\textbf{(III)} Given $T_0$ and $R$, it remains to show the uniform estimate (\ref{n102}) with the ${}_0W_{p,\mu}^{1/2-1/2p}$-norm on the left-hand side. Here one cannot use the seminorm (\ref{n3}) on a finite interval $J=(0,T)$, since the equivalence constants of the norms depend on $T$ as explained above. To overcome this obstacle, observe that (\ref{nl45}) remains valid if one replaces $J$ by the half-line $\R_+$, and let $u,h\in \E_{1,\mu}(J)$ satisfy (\ref{nl43}). We set $u_*:= S u|_{t=0}\in \E_{1,\mu}(\R_+)$, where $S$ is the right-inverse of the temporal trace $\tr_0$ on $\E_{1,\mu}(\R_+)$, and define the functions
$$\wt{u}:= \calE_J^0(u-u_*) + u_* \in \E_{1,\mu}(\R_+), \qquad \wt{h}:= \calE_J^0h\in {}_0\E_{1,\mu}(\R_+),$$
where  $\calE_J^0$ is the extension operator on ${}_0\E_{1,\mu}(J)$ to the half-line, whose operator norm is independent of $T$. Observe that 
$$
|\wt{u}|_{BC([0,\infty)\times \ol{\Omega})} \lesssim |\wt{u}|_{\E_{1,\mu}(\R_+)}\lesssim R + \big |u|_{t=0}\big |_{B_{p,p}^{2(\mu-1/p)}(\Omega)} \lesssim R, 
$$
and, due to $h|_{t=0} = 0$,
$$|\wt{h}|_{BC([0,\infty)\times \ol{\Omega})} \lesssim |\wt{h}|_{\E_{1,\mu}(\R_+)}\lesssim |h|_{\E_{1,\mu}(J)}\leq R,$$ 
where these estimates are independent of $T$. Thus the images  $\wt{u}(\R_+\times \ol{\Omega})$ and $\wt{h}(\R_+\times \ol{\Omega})$ belong to a compact set in $\R^N$, which only depends on $R$, but not on $T\leq T_0$. Thus, using (\ref{nl45}) on the half-line $J=\R_+$, we estimate
\begin{align*}
|g(&u+h)  -g(u) -g'(u)h |_{{}_0W_{p,\mu}^{1/2-1/2p}(J; L_p(\Gamma))} \lesssim |g(\wt{u}+\wt{h}) -g(\wt{u}) -g'(\wt{u})\wt{h}|_{W_{p,\mu}^{1/2-1/2p}(\R_+; L_p(\Gamma))}\\
& \leq \eps(|\wt{h}|_{BC([0,\infty)\times \Gamma)}) \big ( [\wt{h}]_{W_{p,\mu}^{1/2-1/2p}(\R_+; L_p(\Gamma))} + |\wt{h}|_{BC([0,\infty)\times \Gamma)} [\wt{u}]_{W_{p,\mu}^{1/2-1/2p}(\R_+; L_p(\Gamma))}\big )\\
& \lesssim \eps(|h|_{\E_{1,\mu}(J)}) |h|_{\E_{1,\mu}(J)},
\end{align*}
where $\eps$ is uniform in $T\leq T_0$ and $R$. This shows (\ref{n102}). \eprf

As a consequence we have the following result for the boundary map $B$, given by
$$B(u):= \alpha_{ij} \nu_i \tr_\Omega \partial_j u - \big( a^{-1}g\big)(\tr_\Omega u), \qquad u\in \E_{1,\mu}(J).$$

\begin{lem}\label{nl23}\textsl{Let $J=(0,T)$ be finite, and let $p\in (n+2,\infty)$, $\mu\in (1/p,1]$ satisfy $2(\mu-1/p)> 1+n/p$. Then it holds
$$B\in C^1\big (\E_{1,\mu}(J), \F_{\mu}(J)\big ), \qquad B'(u) = \alpha_{ij} \nu_i \emph{\tr}_\Omega \partial_j - \big(a^{-1}g\big)'(\emph{\tr}_\Omega u)\emph{\tr}_\Omega.$$
Further, let $T_0, R>0$ be given. Then there is a continuous function $\eps:[0,\infty)\ra [0,\infty)$ with $\eps(0)=0$, such that for $T\leq T_0$ it holds
$$|B(u+h)-B(u)-B'(u)h|_{ {}_0\F_\mu(J)} \leq \eps(|h|_{\E_{1,\mu}(J)}) |h|_{\E_{1,\mu}(J)}$$
for all $u,h\in \E_{1,\mu}(J)$ as in (\ref{nl43}) above.}
\end{lem}

The last result in this section is concerned with the map $B$ on function spaces without time dependence.

\begin{lem}\label{nl46}\textsl{Let $p\in (n+2,\infty)$ and $\mu\in (1/p,1]$ satisfy $2(\mu-1/p) > 1+n/p$. Then we have
$$B\in C^1\big (W_p^{2(\mu-1/p)}(\Omega), W_p^{2(\mu-1/p) - 1-1/p}(\Gamma)\big),$$ 
with derivative $B'(u_0) = \alpha_{ij} \nu_i \emph{\tr}_\Omega \partial_j  - \big(a^{-1}g\big)'(\emph{\tr}_\Omega u_0)\emph{\tr}_\Omega$ for  $u_0 \in W_p^{2(\mu-1/p)}(\Omega).$ Further, if (\ref{nl3}) is valid, then for each $u_0$ the map $B'(u_0)$ is surjective with bounded linear right-inverse.}
\end{lem}
\bprf Using the continuous right-inverse $S$ of the temporal trace $\tr_0$ on $\E_{1,\mu}(0,1)$ we may write $B(u_0) = \tr_{0} B(S u_0)$. The continuity of $\tr_{0}: \F_{\mu}(0,1) \ra W_p^{2(\mu-1/p) - 1-1/p}(\Gamma)$, which is due to Theorem 4.2 of \cite{MS11a}, and Lemma \ref{nl23} yield that $B$ is $C^1$, with derivative as asserted. 

Now suppose that (\ref{nl3}) holds true, and take $u_0\in W_p^{2(\mu-1/p)}(\Omega)$. For a continuous right-inverse of $B'(u_0)$ we intend to apply Proposition 2.5.1 of \cite{Mey10}. To verify the conditions on $B'(u_0)$ required there, consider the operators $\calA :=   - \alpha_{ij}\partial_i \partial_j$ and $\calB:= \alpha_{ij} \nu_i \tr_\Omega \partial_j.$ The assumption (\ref{nl3}) and Theorem 4.4 of \cite{Ama90} yield that $(\calA,\calB)$ satisfies the ellipticity conditions (E) and (LS) required for the application of the result in \cite{Mey10}. Since these conditions are independent of the lower order terms, also $(\calA, B'(u_0))$ satisfies (E) and (LS). One can show as in the proof of Lemma \ref{jesaba1} that
$\big ( a^{-1} g\big)'(\tr_\Omega u_0)\in  W_p^{2(\mu-1/p) -1/p}(\Gamma),$ which yields that also the required regularity of the coefficients is satisfied. The existence of a continuous right-inverse follows. 
\eprf

\section{The Local Solution Semiflow} \label{locwp}
Using the nonlinear maps $A$ and $B$ defined above we may rewrite (\ref{nl1}) into the form
\begin{alignat}{3}
u_t + A(u) & =    0 & \qquad & \trm{in } \Omega, & \qquad &  t>0,  \nonumber\\
B(u)   & =  0 &&\trm{on } \Gamma,  &&  t>0,\quad  \label{nl14} \\
u|_{t=0} & =  u_0 &&  \trm{in } \Omega. && \nonumber 
\end{alignat} 
For $p\in (n+2,\infty)$ and $s\in (1+n/p,2-2/p]$ we introduce the nonlinear phase spaces
$$\calM_p^s= \big \{ u_0\in W_p^s(\Omega)\;:\; B(u_0)  = 0\big \},$$
which are equipped with the metric of $W_p^s(\Omega)$. We aim to show that (\ref{nl14}) generates a \textsl{compact local semiflow} of $\E_{1,\mu}$-solutions on all these $\calM_p^s$, where $\mu\in (1/p,1]$ is such that $s = 2(\mu-1/p)$. For this the following must be satisfied.
\begin{enumerate}
 \item For all $u_0\in \calM_p^s$ there is $t^+(u_0)>0$ such that (\ref{nl14}) has a unique maximal solution  $u(\cdot,u_0)\in C\big ([0,t^+(u_0)); \calM_p^s\big )$
which belongs to $\E_{1,\mu}(0,\tau)$ for all $\tau\in (0,t^+(u_0))$. 
\item For all $u_0\in \calM_p^s$ and $\tau\in (0,t^+(u_0))$ there is $r>0$ such that $t^+(v_0) >\tau$ for all $v_0\in B_r(u_0) \cap \calM_p^s$, and the map $u(\tau, \cdot):B_r(u_0) \cap \calM_p^s \ra \calM_p^s$ is continuous. 
\item If for a bounded set $M\subset \calM_p^s$ there is $\tau >0$ such that $t^+(v_0) > \tau$ for all $v_0\in M$, then $u(\tau,M)$ is relatively compact in $\calM_p^s$.
\end{enumerate}
We first consider the linearization of (\ref{nl14}) in some $u\in \E_{1,\mu}(J)$, and prove that it admits maximal regularity in the weighted $L_p$-setting. This is the key to well-posedness of (\ref{nl14}) in $\calM_p^s$.

\begin{lem}\label{nl5}\textsl{Let $J=(0,T)$ be finite and let $p\in (n+2,\infty)$, $\mu\in (1/p,1]$ satisfy $s= 2(\mu-1/p) > 1+n/p.$
Assume that  (\ref{nl3}) is valid, let $u\in \E_{1,\mu}(J)$ be given and set
\begin{align*}
 \calD_{u}(J):= \big\{ (\wt{f}, \wt{g}, \wt{v}_0) \in \E_{0,\mu}&(J)\times \F_\mu(J) \times W_p^{s}(\Omega)\;:\;B'\big(u|_{t=0} \big)\wt{v}_0 = \wt{g}_0\;\; \text{ \emph{on} } \;\Gamma \big \}
\end{align*}
Then the linear inhomogeneous, nonautonomous problem 
\begin{alignat}{3}
v_t + A'\big (u(t,x)\big)v & =    \wt{f}(t,x), &\qquad & x\in \Omega, \qquad  &t\in J,   \nonumber\\
B'\big(u(t,x)\big) v  & = \wt{g}(t,x),&&x\in  \Gamma, & t\in J,  \label{n5} \\
v(0,x) & =  \wt{v}_0(x), &&  x\in \Omega. && \nonumber
\end{alignat}
has a unique solution $v\in \E_{1,\mu}(J)$ if and only if the data satisfies $(\wt{f}, \wt{g}, \wt{v}_0) \in \calD_u(J)$. The corresponding solution operator $\calS: \calD_{u}(J)\ra \E_{1,\mu}(J)$ is continuous; i.e., there is $C>0$ with
$$|v|_{\E_{1,\mu}(J)} \leq C \big ( |\wt{f}|_{\E_{0,\mu}(J)} + |\wt{g}|_{\F_\mu(J)} + |\wt{v}_0|_{W_p^{s}(\Omega)}\big)$$
for all $(\wt{f}, \wt{g}, \wt{v}_0) \in \calD_{u}(J)$. Given $T_0>0$, the norm of $\calS$ restricted to $\calD_{u}^0(J) := \big\{ (\wt{f}, \wt{g}, \wt{v}_0) \in \calD_{u}(J)\;:\; \wt{g}\in {}_0\F_\mu(J)\big\}$ has a uniform bound with respect to all $T\leq T_0$.}
\end{lem}
\bprf We intend to apply Theorem 2.1 of \cite{MS11b} to (\ref{n5}), and therefore have to check that $\big(A'(u), B'(u)\big)$ satisfies the regularity conditions (SD), (SB) and the ellipticity conditions (E), (LS$_{\text{stat}}$)  required there. 

We have $u\in C\big (\ol{J}; C^1(\ol{\Omega})\big )$ by (\ref{nl6}), and thus the top order coefficients of $A'(u)$ belong to $C(\ol{J}\times \ol{\Omega})$ and its lower order coefficients belong to $\E_{0,\mu}(J)$. Lemma \ref{jesaba1} implies that the coefficients of $B'(u)$ belong to $\F_{\mu}(J)$. Since the condition $1/2-1/2p > 1-\mu+1/p + \frac{n-1}{2p}$ is equivalent to $2(\mu-1/p)> 1+n/p$ we obtain that (SD) and (SB) of \cite{MS11b} hold true. To verify the ellipticity conditions, consider the operators $\calA$ and $\calB$, given by 
$$\calA v := -\partial_i \big ( a_{ij}(u)  \partial_jv), \qquad \calB v := \alpha_{ij}\nu_i \tr_\Omega  \partial_jv, \qquad v\in \E_{1,\mu}(J).$$
It is shown in Theorem 4.4 of \cite{Ama90} that (\ref{nl3}) implies (E) and (LS$_{\text{stat}}$) of \cite{MS11b} for  $(\calA, \calB)$. Since these conditions are independent of lower order terms, it follows that $\big(A'(u), B'(u)\big)$ satisfies (E) and (LS$_{\text{stat}}$) as well. \eprf

Now we can prove local existence and uniqueness for solutions of (\ref{nl1}) in the weighted setting. The proof is based on the above linear maximal regularity result and the contraction principle, and follows \cite{Zac06} (see also \cite{KPW10} and  \cite{LPS06}).

\begin{lem}\label{nl2}\textsl{Let $p\in (n+2,\infty)$ and $\mu\in (1/p,1]$ satisfy $s= 2(\mu-1/p)>1+n/p$, and assume that (\ref{nl3}) is valid. Then for each initial value $u_0\in \calM_p^s$ the problem (\ref{nl14}) has a unique maximal solution  $u(\cdot,u_0)\in C\big ([0,t^+(u_0)); \calM_p^s\big )$
which belongs to $\E_{1,\mu}(0,\tau)$ for all $\tau\in (0,t^+(u_0))$. Here $t^+(u_0)>0$ denotes the maximal existence time.
}
\end{lem}
\bprf We fix $u_* = S u_0 \in \E_{1,\mu}(\R_+)$ and consider the linearized problem
\begin{alignat}{3}
w_t + A'(u_*)w & =   A'(u_*) u_* - A(u_*) &\qquad & \trm{in } \Omega, \qquad  &t>0,   \nonumber\\
B'(u_*)w & = B'(u_*)u_* - B(u_*)  &&\trm{on } \Gamma, & t> 0,  \label{eq:jezaba9} \\
w|_{t=0} & =  u_0 &&  \trm{in } \Omega. && \nonumber
\end{alignat}
Due to the Lemmas \ref{nl19} and \ref{nl23} it holds $A'(u_*) u_* - A(u_*) \in \E_{0,\mu}(0,1)$ and $B'(u_*)u_* - B(u_*) \in \F_\mu(0,1),$
and since $B(u_0)= 0$ the compatibility condition
$B'(u_0)u_0  =  B'(u_0)u_0 - B(u_0)$ on $\Gamma$ is trivially satisfied. Thus Lemma \ref{nl5} yields a unique solution $w_*\in \E_{1,\mu}(0,1)$ of (\ref{eq:jezaba9}). Using $w_*$, we consider for $\sigma, \tau\in (0,1]$ the closed space
$$\Sigma(\sigma, \tau) := \big \{u\in \E_{1,\mu}(0,\tau)\;:\; |u-w_*|_{\E_{1,\mu}(0,\tau)} \leq \sigma, \;\; u|_{t=0} = u_0\big \}.$$ 
It then follows from the embedding (\ref{nl6}) that
\beq\label{nl24}
|u|_{C([0,\tau]; C^1(\ol{\Omega}))}, \big |u|_{t=0}\big |_{W_p^{2(\mu-1/p)}(\Omega)}, |u|_{\E_{1,\mu}(0,\tau)} \lesssim 1+ |w_*|_{\E_{1,\mu}(0,1)},
\eeq
uniformly in $u \in\Sigma(\sigma,\tau)$ and $\sigma,\tau \in (0,1]$. For $u\in \Sigma(\sigma, \tau)$ we next consider 
\begin{alignat}{3}
w_t + A'(u_*)w & =   A'(u_*)u - A(u) &\qquad & \trm{on } \Omega\times (0,\tau),   \nonumber\\
B'(u_*) w & = B'(u_*)u- B(u)  &&\trm{on } \Gamma\times (0,\tau), \label{eq:jezaba10} \\
w|_{t=0} & =  u_0 &&  \trm{in } \Omega. & \nonumber
\end{alignat}
As above, for all $\tau\in (0,1]$ there is a unique solution $w=\calL(u)\in \E_{1,\mu}(0,\tau)$ of (\ref{eq:jezaba10}) due to Lemma \ref{nl5}. This defines a map $\calL:\Sigma(\sigma,\tau)\ra \E_{1,\mu}(0,\tau).$
Using (\ref{nl24}) and the uniform approximation of the nonlinearities by their derivatives derived in the Lemmas \ref{nl19} and  \ref{nl23}, it is straightforward to see that the contraction principle yields a unique fixed point $u\in \E_{1,\mu}(J)$ of $\calL$ on $\Sigma(\sigma, \tau)$, provided $\sigma$ and $\tau$ are sufficiently small. This fixed point solves (\ref{nl14}). Since for given $\sigma$ each solution of (\ref{nl14}) in $\E_{1,\mu}(0,\tau)$ belongs to $\Sigma(\sigma, \tau)$ for sufficiently small $\tau$, it is in fact the unique solution of (\ref{nl14}). The existence of a maximal existence time $t^+(u_0)$ and a maximal solution in $C\big([0,t^+(u_0)); \calM_p^s\big)$ follows from standard arguments as  e.g. in \cite{LPS06}.\eprf

We next consider the uniformity of local existence times and the continuous dependence of solutions on the initial data. The proof is based on a combination of  maximal regularity  and the implicit function theorem and follows the arguments used in Theorem 14 of \cite{LPS06}.

\begin{lem}\label{nl13}\textsl{In the situation of Lemma \ref{nl2}, let $u = u(\cdot,u_0)$ be the maximal solution of (\ref{nl1}) with initial value $u_0\in \calM_p^s$. Then for all  $\tau\in (0,t^+(u_0))$ there is a ball $B_r(u_0)$ in $W_p^s(\Omega)$, $r>0$, and a continuous map 
$$\Phi:B_r(u_0)\cap \calM_p^s\ra  \E_{1,\mu}(0,\tau), \qquad \Phi(u_0) = u,$$ such that $\Phi(v_0)$ is the solution of (\ref{nl14}) on $(0,\tau)$ with initial value $v_0\in B_r(u_0)\cap \calM_p^s$.}
\end{lem}
\bprf Take $p\in (n+2,\infty)$ and $\mu\in (1/p,1]$ with $s= 2(\mu-1/p)$, such that $u\in \E_{1,\mu}(0,\tau)$. We consider the linear problem
\begin{alignat}{3}
z_t + A'\big (u(t,x)\big)z & =    \wt{f}(t,x), &\qquad & x\in \Omega, \qquad  &t\in (0,\tau),   \nonumber\\
B'\big(u(t,x)\big) z  & = \wt{g}(t,x),&&x\in  \Gamma, & t\in (0,\tau),  \label{jesaba3} \\
z(0,x) & =  \wt{w}_0(x), &&  x\in \Omega, && \nonumber
\end{alignat}
and denote by 
$\calS: \calD_u(0,\tau) \ra \E_{1,\mu}(0,\tau)$ the continuous linear solution operator from Lemma \ref{nl5} corresponding to (\ref{jesaba3}). We have that $v\in \E_{1,\mu}(0,\tau)$ solves (\ref{nl14}) with initial value $v_0\in  \calM_p^s$ if and only if 
\beq\label{nl15}
v = u + \calS\big ( F(v-u), G(v-u), v_0-u_0\big ),
\eeq
where the nonlinear functions $F$ and $G$ are given by
$$F(w) :=  - \big ( A(u+w)- A(u) - A'(u)w\big ), \qquad G(w):=   -\big ( B(u+w) - B(u) - B'(u)w\big ).$$
The Lemmas \ref{nl19} and \ref{nl23} yield that $F\in C^1\big ( \E_{1,\mu}(0,\tau), \E_{0,\mu}(0,\tau) \big )$, $G\in C^1\big ( \E_{1,\mu}(0,\tau), \F_{\mu}(0,\tau) \big ).$ The tangential space of $\calM_p^s$  at $u_0$, which is a closed subspace of $W_p^s(\Omega)$, is given by 
$$T_{u_0}\calM_p^s:= \big \{ z_0 \in W_p^s(\Omega)\;:\; B'(u_0)z_0 = 0\big \}.$$
We consider the nonlinear map $\calF: T_{u_0}\calM_p^s\times \E_{1,\mu}(0,\tau) \ra \E_{1,\mu}(0,\tau),$ defined by 
$$\calF(z_0,w) := w - \calS\big (F(w), G(w), z_0 +  \calN_s \tr_{0}G(w)\big ).$$ 
Here $\calN_s: W_p^{s-1-1/p}(\Gamma) \ra  W_p^s(\Omega)$  denotes the continuous right-inverse of  $B'(u_0)$, which is given by Lemma \ref{nl46}. The map $\calF$ is well defined, since due to $B'(u_0)\big ( z_0 +  \calN_s \tr_{0} G(w)\big )  = \tr_{0}G(w)$
only compatible data are inserted into $\calS$. It further holds $\calF(0,0) = 0$ and that $\calF$ is continuously differentiable. The derivative of $\calF$  with respect to the second argument at $(z_0, w) = (0,0)$ is given by
\begin{align*}
\id + \calS\big ( A'(u+w) - A'(u), B'(u+w) - B'(u), \calN_s \tr_{0}( B'(u+w) - B'(u))\big )|_{w=0}= \id,
\end{align*}
and is therefore invertible. Thus we can  solve the nonlinear equation $\calF(z_0,w)= 0$ locally around $(0,0)$ uniquely by $w= \Phi_*(z_0)$ with a $C^1$-function $\Phi_*:B_r(0) \ra \E_{1,\mu}(0,\tau)$, where  $B_r(0)\subset T_{u_0}\calM_p^s$ with sufficiently small $r>0$.

\textbf{(III)} Now let $v_0\in \calM_p^s$ be given, and define $z_0 := \big ( \id- \calN_sB'(u_0)\big )(v_0-u_0)\in T_{u_0}\calM_p^s.$
By the continuity of $\id- \calN_sB'(u_0)$, if $v_0$ is close to $u_0$ in $\calM_p^s$ then the norm of $z_0$ in $W_p^s(\Omega)$ is small, such that $w = \Phi_* (z_0)\in \E_{1,\mu}(J)$ is well-defined and satisfies
$$w = \calS\big (F(w), G(w), v_0-u_0 - \calN_s (B'(u_0)(v_0-u_0)- \tr_{0}G(w))\big ).$$
Due to $\tr_{0}G(w) = - B(u_0+\tr_{0} w) + B'(u_0)( \tr_{0}w)$, the continuity of $\calN_s$, $B(v_0) = 0$ and  Lemma \ref{nl46} yield
\begin{align}
|\tr_0w - (v_0-u_0)&|_{W_p^s(\Omega)}  = | \calN_s \big (B(u_0+\tr_0w) - B'(u_0)(\tr_0w - (v_0-u_0))\big )|_{W_p^s(\Omega)}\nonumber\\
&\,\quad  \lesssim |B(u_0+\tr_0w) - B(v_0) - B'(v_0)(\tr_0w-(v_0-u_0))|_{W_p^{s-1-1/p}(\Omega)}\nonumber\\
&\,\quad \qquad +  | \big( B'(v_0)- B'(u_0)\big )(\tr_0w-(v_0-u_0))|_{W_p^{s-1-1/p}(\Omega)} \nonumber\\
&\, \quad  \leq \eps\big(|\tr_0w-(v_0-u_0)|_{W_p^s(\Omega)} + |v_0-u_0|_{W_p^s(\Omega)}\big)|\tr_0w-(v_0-u_0)|_{W_p^s(\Omega)},\nonumber
\end{align}
where $\eps$ is a nonnegative continuous function with $\eps(0)=0$. Since $\Phi_*$ is continuous and satisfies $\Phi_*(0) = 0$, if $v_0$ tends to $u_0$ then $|\tr_0w|_{W_p^s(\Omega)}$ tends to zero. Thus for $v_0$ sufficiently close to $u_0$ the above inequality is only possible if $\tr_0 w = v_0-u_0$.
This implies that the function $v= u+w\in \E_{1,\mu}(0,\tau)$ solves (\ref{nl15}), and therefore (\ref{nl14}) with initial value $v_0$. Now 
$$\Phi(v_0):= u + \Phi_*\big ( (\id- \calN_sB'(u_0))(v_0-u_0) \big )$$ 
is the asserted continuous solution map for (\ref{nl14}) on $B_r(u_0)\cap \calM_p^s$. \eprf

The proof shows that the solution map $\Phi$ enjoys in fact more regularity. This can be useful for determining the fractal dimension of an attractor, see e.g. Theorem 4.3 of  \cite{Lad91}.

The above lemma yields that (\ref{nl14}) also  satisfies the second condition required for a local semiflow. We finally show compactness of the flow, employing the inherent smoothing effect of the $L_{p,\mu}$-spaces. Our arguments are inspired by those of Section 3 in  \cite{KPW10}.

\begin{lem}\label{nl50}\textsl{In the setting of Lemma \ref{nl2}, let the bounded set $M\subset \calM_p^s$ and $\tau>0$ satisfy $t^+(v_0)>\tau$ for all $v_0\in M$. Then $u(\tau, M)$ is relatively compact in $\calM_p^s$.}
\end{lem}
\bprf Since the embedding $W_p^s(\Omega) \hra W_p^{s_*}(\Omega)$ is compact for $s_*\in (0, s)$ the set $M$ is relatively compact in $W_p^{s_*}(\Omega)$. Take $\mu_*\in (1/p,1]$ with $s_*= 2(\mu_*-1/p)> 1+n/p$. Due to Lemma \ref{nl13}, for each $v_0\in M$ there is a ball $B_{r}(v_0)$ in $W_p^{s_*}(\Omega)$ and a continuous map $\Phi: B_{r}(v_0) \cap \calM_p^{s_*} \ra \E_{1,\mu_*}(0,\tau)$	
such that $w = \Phi(w_0)  \in \E_{1,\mu_*}(0,\tau)$ solves (\ref{nl14}) with initial value $w_0\in B_{r}(v_0)\cap \calM_p^{s_*}$. This yields an open cover of $M$ in $W_p^{s_*}(\Omega)$, and thus, by compactness, there are finitely many balls $B_{k}$ and maps $\Phi_{k}$ with the above property such that $\bigcup_k B_k$ covers $M$. Each $\Phi_{k}$ maps the relatively compact set $B_{k}\cap M$ continuously into $\E_{1,\mu_*}(0,\tau)$, with $\Phi_{k}(w_0) = u(\cdot,w_0)|_{(0,\tau)}$ for $w_0\in B_{k}\cap M.$ Since the temporal trace $\tr_{\tau}:\E_{1,\mu_*}(0,\tau) \ra W_p^{2-2/p}(\Omega)$, i.e., $\tr_{\tau}w = w|_{t=\tau},$ is continuous, we obtain that $u(\tau,M) =  \bigcup_k \tr_{\tau} \circ \Phi_{k}( B_{k}\cap M)$
is relatively compact in $W_p^s(\Omega),$ as a continuous image of a relatively compact set. 
\eprf

We summarize the above considerations to the main result of this section.

\begin{prop}\label{nl32}\textsl{Let $p\in (n+2,\infty)$, $\mu\in (1/p,1]$ and $s\in (1+n/p,2-2/p]$ satisfy $s= 2(\mu-1/p)$, and assume that (\ref{nl3}) holds true. Then the system (\ref{nl1}) generates a compact local semiflow  of $\E_{1,\mu}$-solutions on the phase space $\calM_p^s$.}
\end{prop}

\begin{rem}\label{remark}\emph{The methods in this section are independent of the concrete form of the nonlinear operators $A$ and $B$, as long as they are $C^1$ and Theorem 2.1 of \cite{MS11b} and Proposition 2.5.1 of \cite{Mey10} are applicable to the corresponding linearization. Thus a compact local semiflow in a scale of nonlinear phase spaces can be obtain for more general, also higher order parabolic systems with nonlinear boundary conditions, as treated in \cite{LPS06}, for instance.}
\end{rem}

\section{Attractors in Stronger Norms} \label{gaisn}
We now assume the situation of Proposition \ref{nl32}, fix $p\in (n+2,\infty)$ and investigate the long-time behaviour of the solution semiflow generated by (\ref{nl1}) in $\calM_p^{2-2/p}$ in terms of attractors. 

Assuming that all solutions of (\ref{nl1}) are global in time, a subset $\calA$ of $\calM_p^{2-2/p}$ is called a \textsl{global attractor} for (\ref{nl1}) if $\calA$ is nonempty, compact, invariant with respect to the semiflow and attracts every bounded subset of $\calM_p^{2-2/p}$. The latter means that for every bounded set $B\subset \calM_p^{2-2/p}$ it holds that
$$d(u(t,B), \calA) = \sup_{v\in u(t,B)} \inf_{w\in \calA} \|v-w\|_{W_p^{2-2/p}(\Omega)} \to 0 \qquad \text{as }\, t\to +\infty,$$
where $u(t,\cdot)$ it the solution operator for  (\ref{nl1}) at time $t>0$.

Using the full strength of maximal $L_{p,\mu}$-regularity we can estimate solutions of (\ref{nl1}) at a later time in a strong norm by the solution at an earlier time in a weaker norm. This estimate is the key to global attractors in stronger norms.

\begin{lem}\label{nl31}\textsl{Let $u_0\in \calM_p^{2-2/p}$, denote by $u(\cdot,u_0)$ the maximal solution of (\ref{nl1}) and let $q\in (1,p]$, $\sigma \in (0,2-2/q]$. Take further $\tau>0$ and $0<T_1<T_2<t^+(u_0)$ with $\tau= T_2-T_1$. Then for $\alpha>0$ there is a constant $C = C\big(\tau,\alpha, |u(\cdot,u_0)|_{C([T_1,T_2], C^\alpha(\ol{\Omega}))}\big)$ such that
\begin{align}
|u(T_2,u_0)|_{W_q^{2-2/q}(\Omega)} \leq C \big ( 1 + |u(T_1,u_0)|_{W_q^{\sigma}(\Omega)} \big ).\label{nl28}
\end{align}
In the semilinear case, i.e., if $(a_{ij})$ does not depend on $u$, one may take $\alpha = 0$.
}
\end{lem}
Let us briefly consider the above estimate in more detail. The main point is that $q$ may be arbitrarily large and that $\sigma$ may be arbitrarily small. Hence, given numbers $\alpha,\beta\in (0,1)$, it follows from (\ref{nl28}) and Sobolev's embeddings that
$$|u(T_2,u_0)|_{C^{1+\beta}(\ol{\Omega})} \leq  C\big(\tau,\alpha, |u(\cdot,u_0)|_{C([T_1,T_2]; C^\alpha(\ol{\Omega}))}\big).$$
We can therefore control the spatial gradient of the solution in a H\"older norm by the solution itself. This is why we call (\ref{nl28}) a \textsl{gradient estimate}. Usually estimates of this type are obtained for small $\beta$; cf. \cite{DiB93}.  Here $\beta$ may be close to $1$. We would finally like to emphasize that the systems under consideration only have to satisfy the general assumptions (\ref{nl3}). \vspace{0.1cm}\\
\noindent\textbf{Proof of Lemma \ref{nl31}.} Throughout we set $J:= (0,\tau)$ and take $\mu\in (1/q,1]$ with $\sigma = 2(\mu-1/q)$. The spaces $\E_{1,\mu}$, $\E_{0,\mu}$  and $\F_\mu$ must now be understood with respect to $q$; e.g., $\E_{0,\mu}(J) = L_{q,\mu}\big (J; L_q(\Omega)\big)$.

\textbf{(I)} Define the function $v\in \E_{1,1}(J)$ by $v(t,\cdot) := u(t + T_1, u_0)$ for $t\in J$. Since the weight only has an effect at $t=0$, the continuity of the trace at $t=\tau$ (Theorem 4.2 of \cite{MS11a}) yields
\beq\label{nl30}
|u(T_2,u_0)|_{W_q^{2-2/q}(\Omega)} = \big |v|_{t=\tau}\big |_{W_q^{2-2/q}(\Omega)} \lesssim |v|_{\E_{1,\mu}(J)}.
\eeq
Moreover, the function $v$ solves the nonautonomous, inhomogeneous linear problem
\begin{alignat}{3}
w_t - a_{ij}(v)\partial_i\partial_jw & =  a_{ij}'(v)\partial_i v \partial_j v +  f(v) & \qquad & \trm{in } \Omega, & \qquad &  t\in J,  \nonumber\\
\alpha_{ij}\nu_i \partial_j w   & =  a^{-1}g(v) &&\trm{on } \Gamma,  &&  t\in J,\quad  \nonumber \\
w|_{t=0} & =  u(T_1,u_0) &&  \trm{in } \Omega. && \nonumber 
\end{alignat} 
As in the proof of Lemma \ref{nl5}, we infer from Theorem 2.1 of \cite{MS11b} that this problem enjoys maximal regularity in the space $\E_{1,\mu}(J)$. A compactness argument thus yields that there is a constant $C$, which is for given $R>0$ uniform for all $u$ such that  $|u|_{C([T_1,T_2]\times \ol{\Omega})}\leq R$, such that
\beq\label{nl29}
|v|_{\E_{1,\mu}(J)}\leq C\big (|a_{ij}'(v)\partial_i v \partial_j v|_{\E_{0,\mu}(J)} +  |f(v)|_{\E_{0,\mu}(J)} + |a^{-1}g(\tr_\Omega v)|_{\F_\mu(J)} + |u(T_1,u_0)|_{W_q^{\sigma}(\Omega)}\big).
\eeq
\textbf{(II)} Using H\"older's inequality, we estimate for the first summand in (\ref{nl29}) 
\begin{align*}
|a_{ij}'(v)\partial_i v \partial_j v|_{\E_{0,\mu}(J)}^q &\, \leq C |\partial_i v \partial_j v|_{\E_{0,\mu}(J)}^q\\
&\,\leq C \big | |\partial_i v|_{L_{2q}(\Omega)} |\partial_j v|_{L_{2q}(\Omega)} \big |_{L_{q,\mu}(J)}^q\leq C \int_J t^{q(1-\mu)} |v(t,\cdot)|_{W_{2q}^1(\Omega)}^{2q} \D t,
\end{align*}
where $C$ is as above. By the fractional order Gagliardo-Nirenberg inequality, shown in Proposition 4.1 of \cite{Ama85}, we have for $t\in J$ that
$$|v(t,\cdot)|_{W_{2q}^1(\Omega)}^{2q} \lesssim |v(t,\cdot)|_{W_{q}^{\vartheta}(\Omega)}^q \,|v(t,\cdot)|_{W_{r}^{\tau}(\Omega)}^q$$
for $r\in (1,\infty)$ and  $\vartheta, \tau >0$, provided $1-\frac{n}{2q} < \frac{1}{2}\big (\tau- \frac{n}{r}\big )+\frac{1}{2}\big (\vartheta- \frac{n}{q}\big ).$ For given $\alpha$ it holds $C^\alpha(\ol{\Omega}) \hra W_r^{\tau}(\Omega)$ for $\tau\in (0,\alpha)$ and $r\in (1,\infty)$. Thus if $\vartheta<2$ is sufficiently close to $2$ and $r$ is large we obtain 
$$|v(t,\cdot)|_{W_{2q}^1(\Omega)}^{2q}  \lesssim |v(t,\cdot)|_{W_{q}^{\vartheta}(\Omega)}^q \,|v(t,\cdot)|_{C^\alpha(\ol{\Omega}) }^q.$$
We now use that  $W_{q}^{\vartheta}$ may be represented as a real interpolation space between $L_q$ and $W_q^2$ for $\vartheta \neq 1$; i.e., $W_q^\vartheta = (L_q, W_q^2)_{\vartheta/2,q}$  (see Theorem 4.3.1/1 of \cite{Tri94}). The interpolation inequality (Theorem 1.3.3 of \cite{Tri94}) thus yields
$|v(t,\cdot)|_{W_{q}^{\vartheta}(\Omega)} \lesssim |v(t,\cdot)|_{W_{q}^{2}(\Omega)}^{\vartheta/2}|v(t,\cdot)|_{L_q(\Omega)}^{1-\vartheta/2}.$
Applying Young's inequality and $C(\ol{\Omega}) \hra L_q(\Omega)$, it follows that
\begin{align*}
|v(t,\cdot)|_{W_{2q}^1(\Omega)}^{2q} \leq  \eps |v(t,\cdot)|_{W_{q}^{2}(\Omega)}^q +  C(\eps, |u|_{C([T_1,T_2]; C^\alpha(\ol{\Omega}))}),
\end{align*}
where $\eps>0$ may be chosen arbitrary small. We therefore have
$$|a_{ij}'(v)\partial_i v \partial_j v|_{\E_{0,\mu}(J)} \leq  \eps |v|_{\E_{1,\mu}(J)} + C(\tau, \eps, |u|_{C([T_1,T_2]; C^\alpha(\ol{\Omega}))})$$
for the first summand in (\ref{nl29}). Note that this term does not occur in the semilinear case.

\textbf{(III)} For the second summand in (\ref{nl29})  it is clear that $|f(v)|_{\E_{0,\mu}(J)}  \leq C(\tau) |f(u)|_{C([T_1,T_2]\times \ol{\Omega})}.$
For the third summand, Lemma \ref{jesaba1} and the mapping properties of the spatial trace $\tr_\Omega$ (Theorem 4.5 of \cite{MS11a}) yield
\begin{align*}
 |a^{-1}g(\tr_\Omega v)|_{\F_\mu(J)} &\, \leq C(|u|_{C([T_1,T_2]\times \ol{\Omega})})(  1+ |\tr_\Omega v|_{\F_\mu(J)}) \\
&\, \leq C(|u|_{C([T_1,T_2]\times \ol{\Omega})})(1+ | v|_{H_{q,\mu}^{1/2}(J; L_q(\Omega))\cap L_{q,\mu}(J; W_q^1(\Omega))}).
\end{align*}
Here $H_{q,\mu}^{s}$ is for $s\in (0,1)$ defined by complex interpolation; i.e., $H_{q,\mu}^{s} = [L_{q,\mu}, W_{q,\mu}^1]_{s}$ (see \cite{MS11a}). The interpolation inequality in the complex case (Theorem 1.9.3 in \cite{Tri94}) shows that 
$$|v|_{H_{q,\mu}^{1/2}(J; L_q(\Omega))} \lesssim |v|_{W_{q,\mu}^1(J; L_q(\Omega))}^{1/2} |v|_{\E_{0,\mu}(J)}^{1/2}.$$ 
Moreover, since $L_{q,\mu}(J; W_q^1(\Omega)) = [\E_{0,\mu}(J), L_{q,\mu}(J; W_q^2(\Omega))]_{1/2}$ by Theorem 1.18.4 of \cite{Tri94}, we have
$|v|_{L_{q,\mu}(J; W_q^1(\Omega))} \lesssim  |v|_{L_{q,\mu}(J; W_q^2(\Omega))}^{1/2}|v|_{\E_{0,\mu}(J)}^{1/2}.$ We therefore obtain from Young's inequality that
\begin{align*}
 |a^{-1}g(\tr_\Omega v)|_{\F_\mu(J)} \leq  \eps |v|_{\E_{1,\mu}(J)} + C(\eps, |u|_{C([T_1,T_2]\times \ol{\Omega})}),
\end{align*}
where $\eps$ is arbitrary. If we combine the above estimates with (\ref{nl29}) and choose $\eps$ sufficiently small, then we may subtract $\eps |v|_{\E_{1,\mu}(J)}$ on both sides of the resulting inequality, to obtain
$$|v|_{\E_{1,\mu}(J)}\leq C(\tau, \alpha, |u|_{C([T_1,T_2]; C^\alpha(\ol{\Omega}))}) (1+ |u(T_1,u_0)|_{W_q^{\sigma}(\Omega)}).$$
Together with (\ref{nl30}), this yields the asserted estimate. In the semilinear case the constant does not depend on the H\"older norm of the solution, since then only the terms $|f(v)|_{\E_{0,\mu}(J)}$ and $|a^{-1}g(\tr_\Omega v)|_{\F_\mu(J)}$ in (\ref{nl29}) are estimated. \hfill\rule{1.5ex}{1.5ex}\vspace{0.3cm}

We can now prove the assertion of Theorem \ref{thm} on a global attractor in the quasilinear case.

\begin{prop}\label{nl37}\textsl{Suppose that there are  $\alpha,R>0$ such that for each solution $u(\cdot,u_0)$ of (\ref{nl1}) with initial value $u_0\in \calM_p^{2-2/p}$ it holds $\limsup_{t\ra t^+(u_0)} |u(t,u_0)|_{C^\alpha(\ol{\Omega})} \leq R.$ 
Then (\ref{nl1}) has a global attractor in  $\calM_p^{2-2/p}$.}
\end{prop}
\bprf We first show that $t^+(u_0)= +\infty$ for all $u_0\in \calM_p^{2-2/p}$. Assume the contrary; i.e., $t^+(u_0)<+\infty$. Then Lemma \ref{nl31} and the embedding $C^\alpha(\ol{\Omega})\hra W_p^{\sigma}(\Omega)$ for $\sigma\in (0,\alpha)$ yield
$$\sup_{t\in [0,t^+(u_0))} |u(t,u_0)|_{W_p^{2-2/p}(\Omega)} \leq C(R),$$
which means that the orbit $\{u(t,u_0)\}_{t\in [0,t^+(u_0))}$ is bounded in $W_p^{2-2/p}(\Omega)$. Using that the latter space  embeds compactly into $W_p^{s}(\Omega)$ for $s\in (0,2-2/p)$, we can argue literally as in proof of Theorem 3.1 in \cite{KPW10} to obtain a contradiction to the maximal existence time. Hence $t^+(u_0) = +\infty$. Another application of Lemma \ref{nl31} implies that there is $R_0>0$ with $\limsup_{t\ra +\infty} |u(t,u_0)|_{W_p^{2-2/p}(\Omega)} \leq R_0$
for all initial values $u_0$. Therefore the compact global semiflow generated by (\ref{nl1}) in $\calM_p^{2-2/p}$ has an absorbant ball, and the existence of a global attractor follows from e.g. \cite[Corollary 1.1.6]{CD00}.
\eprf

We next consider the semilinear case with nonlinear boundary conditions. This completes the proof of Theorem \ref{thm}.

\begin{cor}\label{nl206}\textsl{Assume that $(a_{ij})$ does not depend on $u$, and suppose that there are $q\in (1,\infty)$, $\sigma\in (0,2-2/q]$ and  $R>0$ such that for each solution $u(\cdot,u_0)$ of (\ref{nl1}) with  $u_0\in \calM_p^{2-2/p}$ it holds
$\limsup_{t\ra t^+(u_0)} |u(t,u_0)|_{W_q^\sigma(\Omega) \cap C(\ol{\Omega})} \leq R.$ Then (\ref{nl1}) has a global attractor in $\calM_p^{2-2/p}$.}
\end{cor}
\bprf Lemma \ref{nl31} yields a constant $R_0$ with
\beq\label{nl33}
\limsup_{t\ra t^+(u_0)} |u(t,u_0)|_{W_q^{2-2/q}(\Omega)} \leq R_0
\eeq
for all $u_0\in \calM_p^{2-2/p}$. We employ a bootstrapping procedure to show that (\ref{nl33}) remains true if one replaces $W_q^{2-2/q}(\Omega)$ by $C^\alpha(\ol{\Omega})$ with some $\alpha>0$, and $R_0$ by a possibly larger constant. It then follows from Proposition \ref{nl37} that (\ref{nl1}) has a global attractor in $\calM_p^{2-2/p}$ as asserted. If $q > n/2+1$ then Sobolev's embedding yields 
$W_q^{2-2/q}(\Omega) \hra C^\alpha(\ol{\Omega})$ 
for some $\alpha>0$, and we are done in this case. Otherwise, in case $q\in (1, n/2+1)$, we employ $W_q^{2-2/q}(\Omega) \hra W_{q_1}^\tau(\Omega),$
which is valid for some  $\tau>0$ if $q_1\in \big (q, \frac{nq}{n+2-2q}\big )$. Note here that $\frac{nq}{n+2-2q}> q$ for all $n$ and $q \in  (1, n/2+1)$. Another application of Lemma  \ref{nl31} yields (\ref{nl33}) with $W_q^{2-2/q}(\Omega)$ replaced by $W_{q_1}^{2-2/q_1}(\Omega).$ Iteratively, this yields a strictly increasing sequence of numbers $q_k$ as long as $q_k<n/2+1$. But since $q_k \geq \big (\frac{n(1-\delta)}{n +2 - 2q}\big )^k q$ for small $\delta>0$ as long as $q_k<n/2+1$ and $\frac{n}{n +2 - 2q}>1$, the sequence $q_k$ becomes larger than $n/2+1$ after finitely many steps. Thus (\ref{nl33}) holds true with a H\"older norm, and this finishes the proof.\eprf

We now consider special cases of (\ref{nl1}) where uniform a priori estimates of De Giorgi-Nash-Moser type allow a further reduction of the regularity of the absorbant ball. We start with a single equation, $N=1$.

\begin{cor}\textsl{Suppose that for the scalar reaction terms $f,g$ there are $r,C>0$ such that $|f(\zeta)|\leq C |\zeta|^r$ and $g(\zeta)\zeta \leq C(1+|\zeta|^2)$ is valid for all $\zeta \in \R$, and assume that $a\geq \delta >0$. If 
\begin{alignat}{3}
u_t - \emph{\dv}\big( a(u) \nabla u\big) & =    f(u) & \qquad & \emph{\trm{in }} \Omega, & \qquad &  t>0,  \nonumber\\
a(u)\partial_\nu u   & =  g(u) &&\emph{\trm{on }} \Gamma,  &&  t>0,\quad  \label{nl203} \\
u|_{t=0} & =  u_0 &&  \emph{\trm{in }} \Omega, && \nonumber 
\end{alignat} 
admits an absorbant ball in a $L_q(\Omega)$-norm for some $q\geq 1$ with $q> \frac{n}{2}(r-1)$ then it has a global attractor in $\calM_p^{2-2/p}$ for all $p>n+2$.}
\end{cor}
\bprf It is shown in Theorem 1 of \cite{Dun97} that the existence of an absorbant ball in $L_q(\Omega)$ implies an absorbant ball in $C(\ol{\Omega})$.  This in turn yields an absorbant ball in a H\"older norm, which follows, e.g., from  Theorem III.1.3 of \cite{DiB93} or Corollary 4.2 of \cite{Dun00}.\eprf

We next consider for $(u,v)\in \R^2$ quasilinear cross-diffusion systems of the form 
\begin{alignat}{3}
u_t - \dv\big(P(u,v)\nabla u + R(u,v)\nabla v\big) & =    f_1(u,v) & \qquad & \trm{in } \Omega, & \qquad &  t>0,  \nonumber\\
v_t - \dv\big(Q(v)\nabla v\big) & =  f_2(u,v) & \qquad & \trm{in } \Omega, & \qquad &  t>0,  \nonumber\\
\partial_\nu (u,v)   & =  0&&\trm{on } \Gamma,  &&  t>0,\quad  \label{nl202} \\
(u,v)|_{t=0} & =  (u_0,v_0) &&  \trm{in } \Omega. && \nonumber 
\end{alignat} 
This problem fits into our setting with $a(u,v) = \left ( \begin{array}{cc} P(u,v) & R(u,v) \\ 0 & Q(v) \end{array}\right)$, $\alpha_{ij} = \delta_{ij}$ and $g = 0$. We assume that there are nonnegative continuous functions $\Phi_1$, $\Phi_2$ and constants $C,d>0$ such that for all $\zeta=(\zeta_1, \zeta_2)\in \R^2$ it holds
\beq
\left . \begin{array}{c} 
P(\zeta) \geq d (1+\zeta_1), \quad \zeta_1 \geq 0, \qquad |R(\zeta)| \leq \Phi_1(\zeta_2)\zeta_1,  \qquad Q(\zeta_2)\geq d;\\
\text{the partial derivatives of $P$, $R$ are majorized by some powers of $\zeta_1, \zeta_2$;}\\
|f(\zeta)|\leq \Phi_2(\zeta_2) (1+\zeta_1), \quad g(\zeta)\zeta_1^r \leq \Phi_2(\zeta_2) ( 1+ \zeta_1^{r+1}), \quad \text{ for all } \zeta_1,\zeta_2 \geq 0,\; r>0.
\end{array}\right\}\label{n7}
\eeq
With the results of \cite{KD07} we can weaken the regularity for an  absorbing ball of (\ref{nl202}).

\begin{cor}\textsl{Assume that (\ref{n7}) is valid, and let the solutions of (\ref{nl202}) be nonnegative for nonnegative initial data. Suppose that there are $r>n/2$ and $R>0$ such that for all $(u_0,v_0)\in \calM_p^{2-2/p}$ it holds $\limsup_{t\ra t^+(u_0,v_0)} |u(t,u_0)|_{L_r(\Omega)} + |v(t,u_0)|_{C(\ol{\Omega})} \leq R.$ Then (\ref{nl202}) has a global attractor in $\calM_{p,+}^{2-2/p} = \{(u_0,v_0)\in \calM_p^{2-2/p}: u_0,v_0\geq 0\}$. One can take $r=1$ if $Q$ does not depend on $v$ .
}
\end{cor}
\bprf It is shown in Theorems 7 and 8 of \cite{KD07} that for all $p>n+2$ the solution semiflow for (\ref{nl202}) in the phase space $\{(u_0,v_0)\in W_p^1(\Omega,\R^2): u_0,v_0\geq 0\}$ has a global attractor. From this the existence of an absorbant set in $\{(u_0,v_0)\in C^\alpha(\ol{\Omega},\R^2): u_0,v_0\geq 0\}$ for some $\alpha >0$ follows. Arguing as in the proof of Proposition  \ref{nl37} we obtain an attractor in $\calM_{p,+}^{2-2/p}$. \eprf

\section{Applications}\label{app}
We apply the results of the last section to show convergence to attractors in stronger norms for concrete models. Our first example is concerned with semilinear systems of the form
\begin{alignat}{3}
u_t -\Delta u  & =    f(u) & \qquad & \trm{in } \Omega, & \qquad &  t>0,  \nonumber\\
\partial_\nu u   & =  g(u) &&\trm{on } \Gamma,  &&  t>0,\quad  \label{nl35} \\
u|_{t=0} & =  u_0 &&  \trm{in } \Omega, && \nonumber 
\end{alignat}
as considered in \cite{COPR97}. Here the smooth nonlinearities $f,g:\R^N\ra \R^N$ are dissipative in the sense that there are real numbers $c_i$, $d_i$ with
$\limsup_{|\xi_i| \ra \infty} \frac{f_i(\xi)}{\xi_i} < c_i$ and $\limsup_{|\xi_i| \ra \infty} \frac{g_i(\xi)}{\xi_i} < d_i$, $i=1,...,N$,
such that the first eigenvalue $\lambda_0$ of the linear elliptic problem
\begin{alignat*}{3}
-\Delta v  - c v& =    \lambda v & \qquad & \trm{in } \Omega, & \qquad &   \nonumber\\
\partial_\nu v  -dv & =  0 &&\trm{on } \Gamma,  &&  \quad  \nonumber
\end{alignat*}
is positive, where $c = (c_1,...,c_N)$ and $d = (d_1,...,d_N)$. The discussion in Section 6 of \cite{COPR97} shows that the first eigenvalue of the above problem can be positive although $c_i$ or $d_i$ has the `wrong` sign, i.e., is positive. In this sense $f$ can compensate a possible nondissipativeness of $g$, and vice versa. In Theorem 4.1 of \cite{COPR97} it is shown that under the above assumptions (\ref{nl35}) has a global attractor in the linear phase space $W_2^1(\Omega,\R^N)\cap C(\ol{\Omega},\R^N)$. The detailed balance between $f$ and $g$ and more refined conditions for the existence of a global attractor are discussed in \cite{RT01}. Corollary \ref{nl206} improves the result of \cite{COPR97} as follows.
\begin{thm}\textsl{Under the above assumptions, for $p\in (n+2,\infty)$ the semiflow generated by (\ref{nl35}) has a global attractor in the nonlinear phase space
$$\big\{ u_0\in W_p^{2-2/p}(\Omega,\R^N)\;:\; \partial_\nu u_0 = g(u_0)\;\;\;\text{on }\;\Gamma\big \}.$$}
\end{thm}

As a next example we consider a chemotaxis model with volume-filling effect, 
\begin{alignat}{3}
u_t - d_1\Delta u - \dv\big (  uq(u) \chi(v) \nabla  v\big ) & =    uf(u) & \qquad & \trm{in } \Omega, & \qquad &  t>0,  \nonumber\\
v_t - d_2\Delta v & =  g_1(u) - vg_2(v) & \qquad & \trm{in } \Omega, & \qquad &  t>0,  \nonumber\\
\partial_\nu (u,v)  & =  0 &&\trm{on } \Gamma,  &&  t>0,\quad   \nonumber \\
(u, v)|_{t=0} & = (u_0,v_0) &&  \trm{in } \Omega, && \label{nl38}
\end{alignat}
where $(u,v)\in \R^2$. This model has been introduced in \cite{HP01}, and may be cast in the form (\ref{nl202}) as above. It is assumed that $q(u) = 1- u/U_M$ for some $U_M>0$, that $d_1,d_2 > 0$  and 
$$f|_{(U_M,\infty)} \leq 0, \qquad g_1,g_2\geq 0, \qquad g_1(0) = 0, \qquad \lim_{v\ra \infty} vg_2(v) \ra +\infty.$$ 
 In \cite{Wrz04} it is shown that (\ref{nl38}) has under these assumptions a global attractor in the phase spaces
$$\big \{(u_0,v_0)\in W_p^1(\Omega,\R^2)\;:\; 0\leq u_0 \leq U_M,\;\;\; 0\leq v_0 \big \}, \qquad p\in (n,\infty).$$ 
In \cite{JZ09} it is shown that in fact every solution of (\ref{nl38}) converges to an equilibrium. The proof of Proposition \ref{nl37} yields the following improvement of the result in \cite{Wrz04}.
\begin{thm}\textsl{Under the above assumptions, for $p\in (n+2,\infty)$ the chemotaxis model (\ref{nl38}) has  a global attractor in the phase space
$$\big \{(u_0,v_0)\in W_p^{2-2/p}(\Omega,\R^2)\;:\; 0\leq u_0 \leq U_M,\;\;\; 0\leq v_0 \big \}.$$}
\end{thm}

Our last example is the Shigesada-Kawasaki-Teramoto cross-diffusion model for population dynamics, introduced in \cite{SKT79}, which is for $(u,v)\in \R^2$ given by
\begin{alignat}{3}
u_t - \Delta\big (d_1 + \alpha_{11} u + \alpha_{12} v) u\big) & =   u(a_1-b_1u -c_1 v)& \qquad & \trm{in } \Omega, & \qquad &  t>0,  \nonumber\\
v_t - \Delta\big (d_2 + \alpha_{21} u + \alpha_{22} v) v\big) & =   v(a_2-b_2u -c_2 v)& \qquad & \trm{in } \Omega, & \qquad &  t>0,  \nonumber\\
\partial_\nu (u, v)  &=  0 &&\trm{on } \Gamma,  &&  t>0,\quad   \nonumber \\
(u, v)|_{t=0} &= (u_0,v_0) &&  \trm{in } \Omega. && \label{nl207}
\end{alignat}
Again this model may be cast in the form (\ref{nl202}). Here the constants $a_i,b_i, c_i, d_i$, $i=1,2$, are positive, and the constants $\alpha_{ij}$, $i=1,2$, are nonnegative. In Theorem 2 of \cite{KD07} it is shown that (\ref{nl207}) has a global attractor as a dynamical system in $W_p^1(\Omega,\R^2)$ for $p\in (n,\infty)$, provided $\alpha_{22} = 0$. For $n=2$ this remains true also for $\alpha_{22} > 0$. Proposition \ref{nl37} improves this as follows.
\begin{thm}\textsl{Under the above assumptions, for $p\in (n+2,\infty)$ the population model (\ref{nl207})  has a global attractor in the phase space $W_p^{2-2/p}(\Omega,\R^2)$.}\bigskip
\end{thm}

\end{document}